\documentclass[a4paper,12pt,dvipdfmx]{article}
\usepackage{graphicx}

\usepackage{amsmath,amssymb,cite,comment,bm,color}

\def\beq{\begin{equation}}
\def\eqn#1{\beq\label{#1}}
\def\eeq{\end{equation}}

\font\tfont=cmbx12 scaled\magstep1 %large
\def\veps{\varepsilon}

\def\bb {\begin {eqnarray}}
\def\eqnn#1{\bb\label{#1}}
\def\ee {\end {eqnarray}}

\def\rf#1#2{(\ref{#1}{#2})}

\newcommand{\eqna}[1]{\begin{subequations} \label{#1}
\begin{eqnarray}}
\def\eena{\end{eqnarray}
\end{subequations}}

\def\bbz{\mathbb{Z}}%{Z\!\!\!Z}
\def\bbc{\mathbb{C}}%{I\!\!\!\!C}}

\def\bbr{\mathbb{R}}%{{I\!\!R}}
\def\bbn{\mathbb{N}}%{I\!\!N}

\def\nn{\nonumber}
\def\nt{\noindent}

\def\np{\vfill\eject}
 \def\nd{\end{document}}

\def\trh{{\textstyle{3\over2}}}
\def\fha{{\textstyle{5\over2}}}

\def\sevha{{\textstyle{7\over2}}}
\def\llr{\longrightarrow}
\def\({\left(}
\def\){\right)}

\def\lra{\longrightarrow}

\def\dia{{$\diamondsuit$}}

\def\ha{{\textstyle{1\over2}}}

  \def\tV{{\tilde V}}
\def\tcn{{\tilde{\cal N}}}

\def\r{\rho}
\def\cgc{{\cg^\bbc}}

\def\bbr{{I\!\!R}}
\def\bbn{I\!\!N}
\def\a{\alpha}
\def\b{\beta}

\def\vr{\vert}

\def\l{\lambda}

\def\D{{\Delta}}

\def\L{\Lambda}

\def\riga{-\kern-4pt - }
\def\rigal{-\kern-4pt - \kern-4pt -}
\font\fat=cmsy10 scaled\magstep5

\def\Bbullet{\raise-3pt\hbox{\fat\char"0F}}

\def\ca{{\cal A}}  \def\cc{{\cal C}}
\def\cd{{\cal D}} \def\ce{{\cal E}} \def\cf{{\cal F}}
\def\cg{{\cal G}} \def\ch{{\cal H}} 
 \def\ck{{\cal K}} 
\def\cm{{\cal M}} \def\cn{{\cal N}} 
\def\cp{{\cal P}} \def\cq{{\cal Q}} 
 \def\ct{{\cal T}}

\def\tcd{\widetilde{\cal D}}

\def\ido{intertwining differential operator}
\def\idos{intertwining differential operators}

 \def\ha{{\textstyle{\frac{1}{2}}}}

\def\ca{{\cal A}}
\def\nn{\nonumber}

\def\nt{\noindent}

\def\lra{\longleftrightarrow}
\def\lla{\longleftarrow}

\input epsf.tex
\newcount\figno
\def\fig#1#2#3{
\par\begingroup\parindent=0pt\leftskip=1cm\rightskip=1cm\parindent=0pt
\baselineskip=11pt \global\advance\figno by 1 %\midinsert
\epsfxsize=#3 \centerline{\epsfbox{#2}} \vskip 12pt
#1\par
\endgroup\par}
\def\figlabel#1{\xdef#1{\the\figno}}
\def\encadremath#1{\vbox{\hrule\hbox{\vrule\kern8pt\vbox{\kern8pt
\hbox{$\displaystyle #1$}\kern8pt} \kern8pt\vrule}\hrule}}
\begin{document}

\begin{center}

{\tfont
  Heisenberg Parabolic Subgroups of Exceptional\\[4pt] Noncompact\  $G_{2(2)}$\  and %\\[3pt]
    Invariant Differential  Operators}
%Corresponding Verma Modules}

\vskip 1.5cm

{\bf V.K. Dobrev}

 \vskip 5mm

  Institute for Nuclear Research and Nuclear Energy,\\ Bulgarian
Academy of Sciences,\\ 72 Tsarigradsko Chaussee,  1784 Sofia,
Bulgaria

\end{center}

\vskip 1.5cm

 \centerline{{\bf Abstract}}

 In the present paper we continue the project of systematic
construction of invariant differential operators on the example of
the non-compact  algebra  $G_{2(2)}$. We use both the minimal and the maximal Heisenberg
parabolic subalgebras. We give
 the main  multiplets of indecomposable elementary
representations. This includes the explicit parametrization of the \idos\ between the ERs.
 These are new results applicable in all cases when one would like to use  $G_{2(2)}$ invariant differential operators.

\vskip 0.5cm

\nt 

\vskip 1.5cm

\section{Introduction}

Invariant differential operators   play very important role in the
description of physical symmetries.
In a recent paper \cite{Dobinv} we started the systematic explicit
construction of invariant differential operators. We gave an
explicit description of the building blocks, namely, the parabolic
subgroups and subalgebras from which the necessary representations
are induced. Thus we have set the stage for study of different
non-compact groups. An update of the developments as of 2016 is given in \cite{VKD1}
(see also \cite{Dobparab}).

 In the present paper we  focus on the  exceptional non-compact  algebra  $G_{2(2)}$.
 Below we give the general preliminaries necessary for our approach.
%In the next section we give the general preliminaries necessary for our approach.
In Section 2 we introduce the Lie algebra ~$G_2$ (following \cite{Dobg2}), its real from $G_{2(2)}$, and shortly
the corresponding Lie group.
 In Section 3 we consider the representations induced from the minimal parabolic subalgebra of $G_{2(2)}$.
 In Section 4 we consider the representations induced from the two maximal Heisenberg parabolic subalgebras of $G_{2(2)}$.

\bigskip\bigskip

%\section
{\it Preliminaries}

 Let $G$ be a semisimple non-compact Lie group, and $K$ a
maximal compact subgroup of $G$. Then we have an Iwasawa
decomposition ~$G=KA_0N_0$, where ~$A_0$~ is abelian simply
connected vector subgroup of ~$G$, ~$N_0$~ is a nilpotent simply
connected subgroup of ~$G$~ preserved by the action of ~$A_0$.
Further, let $M_0$ be the centralizer of $A_0$ in $K$. Then the
subgroup ~$P_0 ~=~ M_0 A_0 N_0$~ is a minimal parabolic subgroup of
$G$. A parabolic subgroup ~$P ~=~ M A N$~ is any subgroup of $G$
  which contains a minimal parabolic subgroup.

The importance of the parabolic subgroups comes from the fact that
the representations induced from them generate all (admissible)
irreducible representations of $G$ \cite{Lan,Zhea,KnZu}.
   Actually, it may happen that there are
overlaps which should be taken into account.

Let ~$\nu$~ be a (non-unitary) character of ~$A$, ~$\nu\in\ca^*$,
let ~$\mu$~ fix an irreducible representation ~$D^\mu$~ of ~$M$~ on
a vector space ~$V_\mu\,$.

 We call the induced
representation ~$\chi =$ Ind$^G_{P}(\mu\otimes\nu \otimes 1)$~ an
~{\it elementary representation} of $G$ \cite{DMPPT}.
(In the mathematical literature these representations are called "generalised
principal series representations", cf., e.g., \cite{Knapp}.)
  Their spaces of functions are:
\eqn{fun} \cc_\chi ~=~ \{ \cf \in C^\infty(G,V_\mu) ~ \vr ~ \cf
(gman) ~=~ e^{-\nu(H)} \cdot D^\mu(m^{-1})\, \cf (g) \} \eeq where
~$a= \exp(H)\in A$, ~$H\in\ca\,$, ~$m\in M$, ~$n\in N$. The
representation action is the $left$ regular action: \eqn{lrr}
(\ct^\chi(g)\cf) (g') ~=~ \cf (g^{-1}g') ~, \quad g,g'\in G\ .\eeq

An important ingredient in our considerations are the   Verma modules ~$V^\L$~ over
~$\cg^\bbc$, where ~$\L\in (\ch^\bbc)^*$, ~$\ch^\bbc$ is a Cartan
subalgebra of ~$\cg^\bbc$, the weight ~$\L = \L(\chi)$~ is determined
uniquely from $\chi$ \cite{Har2,Dob}.

Actually, since our ERs will be induced from finite-dimensional
representations of ~$\cm$~ (or their limits)  the Verma modules are
always reducible. Thus, it is more convenient to use ~{\it
generalized Verma modules} ~$\tV^\L$~ such that the role of the
highest/lowest weight vector $v_0$ is taken by the
 space ~$V_\mu\,v_0\,$. For the generalized
Verma modules (GVMs) the reducibility is controlled only by the
value of the conformal weight $d$ (related to $\nu$, see below). Relatedly, for the \idos{} only
the reducibility w.r.t. non-compact roots is essential.

One main ingredient of our approach is as follows. We group the
(reducible) ERs with the same Casimirs in sets called ~{\it
multiplets} \cite{Dobmul,Dob}. The multiplet corresponding to fixed
values of the Casimirs may be depicted as a connected graph, the
vertices of which correspond to the reducible ERs and the lines
between the vertices correspond to intertwining operators. The
explicit parametrization of the multiplets and of their ERs is
important for understanding of the situation.

In fact, the multiplets contain explicitly all the data necessary to
construct the \idos{}. Actually, the data for each \ido{} consists
of the pair ~$(\b,m)$, where $\b$ is a (non-compact) positive root
of ~$\cg^\bbc$, ~$m\in\bbn$, such that the BGG \cite{BGG} Verma
module reducibility condition (for highest weight modules) is
fulfilled: \eqn{bggr} (\L+\r, \b^\vee ) ~=~ m \ , \quad \b^\vee
\equiv 2 \b /(\b,\b) \ .\eeq When \eqref{bggr} holds then the Verma
module with shifted weight ~$V^{\L-m\b}$ (or ~$\tV^{\L-m\b}$ ~ for
GVM and $\b$ non-compact) is embedded in the Verma module ~$V^{\L}$
(or ~$\tV^{\L}$). This embedding is realized by a singular vector
~$v_s$~ determined by a polynomial ~$\cp_{m,\b}(\cg^-)$~ in the
universal enveloping algebra ~$(U(\cg_-))\ v_0\,$, ~$\cg^-$~ is the
subalgebra of ~$\cg^\bbc$ generated by the negative root generators
\cite{Dix}.
More explicitly, \cite{Dob}, ~$v^s_{m,\b} = \cp^m_{\b}\, v_0$ (or ~$v^s_{m,\b} = \cp^m_{\b}\, V_\mu\,v_0$ for GVMs).
  Then there exists \cite{Dob} an \ido{} \eqn{lido}
\cd^m_{\b} ~:~ \cc_{\chi(\L)} ~~\llr ~~ \cc_{\chi(\L-m\b)} \eeq given
explicitly by: \eqn{mido}\cd^m_{\b} ~~=~~ \cp^m_{\b}(\widehat{\cg^-})
\eeq where ~~$\widehat{\cg^-}$~~ denotes the $right$ action on the
functions ~~$\cf$, cf. \eqref{fun}.

\section{The non-compact Lie group and algebra of type $G_{2}$}

%Let $G$ denote for the connected, simply connected split real Lie group of type $G_{2}$.
%The maximal compact subgroup is $K = SU(2) \times SU(2)$.

Let ~~$\cgc = G_2$, with Cartan matrix: ~~$(a_{ij}) = \begin{pmatrix}2 & -3 \cr -1 & 2\cr\end{pmatrix}$, simple roots ~~$\a_1,\a_2$~~ with products:
~~$(\a_2, \a_2) = 3(\a_1, \a_1) = -2(\a_2, \a_1)$. We choose
$(\a_1, \a_1) = 2$, ~~then ~~$(\a_2, \a_2) = 6, ~~(\a_2, \a_1) = -3$. ~(Note that in \cite{Dobg2} we have chosen $\a_1$ as long root, $\a_2$ as short root.)~ 
As we know  $G_2$ is 14--dimensional.
 The positive roots may be chosen as:
  \eqn{posg2} \D^+ = \{
\a_1, ~~\a_2, ~~\a_1 + \a_2, ~~\a_2 + 2\a_1, ~~\a_2 + 3\a_1, ~~2\a_2 +
3\a_1 \} \eeq

  We shall use the orthonormal basis
~~$\veps_1, \veps_2, \veps_3$. In its terms
% the simple roots we may choose:
%\eqn{pig2} \a_1 = \veps_2 - \veps_3, \quad \a_2 =
%\veps_1 - 2\veps_2 + \veps_3 \eeq
%\eqn{pig22} \a_1 = \veps_1 - \veps_2, \quad \a_2 =
%-\veps_1 + 2\veps_2 - \veps_3
% \eeq
 the positive roots are given as:
\eqna{noncp}  &&\a_1 = \veps_1 - \veps_2,   ~\a_3 = \a_1 + \a_2=\veps_2 - \veps_3,  ~\a_{4}= \a_2 + 2\a_1 =\veps_1 - \veps_3 ~   \\ &&
\a_2=-\veps_1 + 2\veps_2 - \veps_3, ~ \a_{5} = \a_2 +
3\a_1= 2\veps_1  -\veps_2  - \veps_3, \\ && \a_{6} = 2\a_2 + 3\a_1= \veps_1 + \veps_2 - 2\veps_3
 \nn\eena
 where for future reference we have introduced also notation for the non-simple roots.
 (Note that in \rf{noncp}a are the  short roots,
 in \rf{noncp}b are the long roots.)

%NB Temporary:
%\eqna{noncpt}  &&\g_1=\a_1,   ~\a_3=\g_{11} = \a_1 + \a_2,  ~\a_4= \g_{12}= \a_2 + 2\a_1 ~   \\ &&
%\g_2 = \a_2, ~ \a_5=\g_{13} = \a_2 +
%3\a_1,  ~ \a_6=\g_{23} = 2\a_2 + 3\a_1\eena

%\eqnn{nons} &&\a_3 \equiv \a_1 + \a_2 = \veps_2- \veps_3,\quad  \a_4 \equiv 2\a_1 +
%\a_2 = \veps_1 - \veps_3, \\
%&&\a_5 \equiv 3\a_1 + 2\a_2 = \veps_1 +
%\veps_2 - 2\veps_3, \quad \a_6\equiv 3\a_1 + \a_2 = 2\veps_1 -
%\veps_2 - \veps_3 \nn\ee

Another way to write the roots in general is ~$\b = (b_1,b_2,b_3)$ under the condition   $b_1+b_2+b_3=0$. Then:

%\eqna{ordr} &&\a_1 = (0,1,-1), ~~  \a_2 + \a_1 =(1,-1,0), ~~ \a_2 + 2\a_1=(1,0.-1) \\
%&&\a_2 = (1,-2,1), ~~ \a_2 +3\a_1 =   (1,1,-2), ~~
 %2\a_2 + 3\a_1 = (2,-1,-1), \eena
  %  where in \rf{ordr}a are the short   positive roots, in
%\rf{ordr}b are the long positive roots.

  \eqnn{ordr} &&\a_1 = (1,-1,0), ~~  \a_2 + \a_1 =(0,1,-1), ~~ \a_2 + 2\a_1=(1,0,-1) \\
&&\a_2 = (-1,2,-1), ~~ \a_2 +3\a_1 =   (2,-1,-1), ~~
2\a_2 + 3\a_1 = (1,1,-2)
  \nn\ee

The dual
roots are: ~~$\a_1^\vee ~~=~~ \a_1$, ~~$\a_2^\vee ~~=~~ \a_2/3$,  ~~$(\a_2
+ \a_1)^\vee ~~=~~ \a_2 + \a_1 ~~=~~ 3\a_2^\vee + \a_1^\vee$, ~~$(\a_2
+ 2\a_1)^\vee ~~=~~ \a_2 + 2\a_1 ~~=~~ 3\a_2^\vee + 2\a_1^\vee$,
~~$(\a_2 + 3\a_1)^\vee ~~=~~ (\a_2 + 3\a_1)/3 ~~=~~ \a_2^\vee +
\a_1^\vee$, ~~$(2\a_2 + 3\a_1)^\vee ~~=~~ (2\a_2 + 3\a_1)/3 ~~=~~
2\a_2^\vee + \a_1^\vee$.

The Weyl group ~$W(\cgc,\ch^\bbc)$~ of ~$G_2$~ is the dihedral group of order 12. This
follows from the fact that ~$(s_1\,s_2)^6=1$, where ~$s_1\,,s_2$~
are the two simple reflections.

The complex Lie algebra $G_2$ has one non-compact real form: $\cg = G_{2(2)}$ which is naturally split.
Its maximal compact subalgebra is $\ck = su(2) \oplus su(2)$, also written as  $\ck = su(2)_{S} \oplus su(2)_{L}$ to emphasize
the relation to the root system (after complexification the first factor contains a short root, the second - a long root).
We remind that  $\cg = G_{2(2)}$ has ~{\it discrete series representations}. Actually, it is ~{\it quaternionic discrete series} since ~$\ck$~ contains as direct summand
(at least one) $su(2)$ subalgebra. The number of discrete series is equal to the ratio ~$\vert W(\cgc,\ch^\bbc)\vert / \vert W(\ck^\bbc,\ch^\bbc)\vert$,
where $\ch$ is a compact Cartan subalgebra of both $\cg$ and $\ck$, ~$W$ are the relevant Weyl groups \cite{Knapp}.
Thus, the number of discrete series in our setting is three. They will be identified below.

The compact Cartan subalgebra ~$\ch$~ of $\cg$ will be chosen (following \cite{Vog}) to coincide with the Cartan subalgebra of $\ck$ and
we may write: ~$ \ch ~=~ u(1)_S \oplus u(1)_L\,$. (One may write as in \cite{Vog}  ~$\ch ~=~   \ct_1 \oplus \ct_2$~ to emphasize the torus nature.)
Accordingly, we choose for the positive root system of ~$\ck^\bbc$~  the roots ~$\a_1+\a_2= (0,1,-1)$, and ~$\a_2 + 3\a_1= (2,-1,-1)$ (which are orthogonal to each other).
The lattice of characters of ~$\ch$~ is ~$\l_\ch ~=~ \ha\mu(\a_1+\a_2) + \ha\nu (\a_2 + 3\a_1)$, where ~$\mu,\nu\in \bbz$.

 The complimentary to ~$\ck$~ space is ~$\cq$~ and it is eight-dimensional.

The Iwasawa decomposition of $\cg$ is:
\eqnn{iwas} \cg ~&=&~ \ck \oplus \ca_0 \oplus \cn_0 \\
&&  \dim\, \ca_0 =2, ~~ \dim\, \cn_0 =6 \nn\ee

The Bruhat decomposition is:
\eqnn{bruh}\cg ~&=&~ \tcn_0 \oplus \cm_0  \oplus \ca_0 \oplus \cn_0 \\
 && \cm_0 = 0, ~~ \dim\,  \tcn_0 =6 \nn\ee

 Accordingly the minimal parabolic of $\cg$  is:
\eqn{g2min} \cp_0 ~=~ \cm_0 \oplus \ca_0 \oplus \cn_0 ~=~  \ca_0 \oplus \cn_0 \eeq

\bigskip

 There are two isomorphic maximal cuspidal parabolic subalgebras of $\cg$ which are of Heisenberg type:
 \eqnn{g2max} &&\cp_k ~=~ \cm_k \oplus \ca_k \oplus \cn_k, \quad k=1,2;  \\
 && \cm_k = sl(2,\bbr)_k, ~~ \dim\, \ca_k =1, ~~ \dim\, \cn_k =5
 \nn\ee
Let us denote by ~$\ct_k$~ the compact Cartan subalgebra of ~$\cm_k$. (Recall that ~$[\cm_j,\ca_k] =0$ for $j\neq k$.) Then ~$\ch_k = \ct_k \oplus \ca_k$~ is a non-compact Cartan subalgebra of $\cg$.
 We choose ~$\ct_1$~ to be generated by the short $\ck$-compact root $\a_1+\a_2$  and ~$\ca_1$~  to be generated by the  long root $\a_2$, ~while
~$\ct_2$~ to be generated by the long $\ck$-compact root $\a_2+3\a_1$  and ~$\ca_2$~  to be generated by the  short  root $\a_1$.

 Equivalently,  the $\cm_1$-compact root of ~$\cgc$~ is ~ $\a_1+\a_2$, while the $\cm_2$-compact root is $\a_2+3\a_1$.
 In each case the remaining five positive roots of ~$\cgc$~ are ~$\cm_k$-noncompact.

To characterize the Verma modules we shall use first the Dynkin labels:
    \eqn{dynk} m_i ~\equiv~ (\L+\r,\a^\vee_i)   , \quad i=1,2, \eeq where   ~$\r$~ is
half the sum of the positive roots of ~$\cg^\bbc$.  Thus, we shall use  :
\eqn{dynka} \chi_\L ~=~ \{ m_1, m_2 \} \eeq

Note that when  both ~$m_i\in\bbn$~ then ~$\chi_\L$~ characterizes the finite-dimensional irreps
of ~$\cgc$~ and its real forms, in particular, $\cg$. Furthermore,  ~$m_k\in\bbn$~
characterizes the finite-dimensional irreps  of the $\cm_k$ subalgebra.

 We shall use also the Harish-Chandra parameters:
 \eqn{dynkhc}
 m_{\b}   ~=~ (\L+\r, \b^\vee )\ ,
   \eeq
   for any positive root $\b$,
and explicitly in terms of the Dynkin labels:
\eqna{labhc} \chi_{HC} ~&=&~ \{\ m_{1},  ~~
 m_{3} = 3m_{2}+m_{1},  ~~ m_{4} = 3m_{2}+2m_{1}
 \\ &&
 m_{2}, ~~ m_{5} = m_{2}+m_{1} , ~~m_{6} = 2m_{2}+m_{1} ,  \}
\eena

%%==========================

\section{Induction from minimal parabolic}

\subsection{Main multiplets}

The main multiplets   are in 1-to-1 correspondence
with the finite-dimensional irreps of ~$G_2 $, i.e., they are
labelled by  the   two positive Dynkin labels    ~$m_i\in\bbn$.
When we induce from the minimal parabolic the main multiplets of $\cg$ are the same
as for the complexified Lie algebra $\cgc$. The latter were considered in \cite{Dobg2} but here
we give a different parametrization.

We take ~$\chi_0 = \chi_{HC}$. It has two embedded Verma modules with
HW ~$\L_{1} = \L_0-m_{1}\a_{1}$, and ~$\L_{2} = \L_0-m_{2}\a_{2}$.
The number of ERs/GVMs   in a  main multiplet   is  $12 = \vert W(\cgc)\vert$.
We give the whole multiplet as follows:
\eqnn{enum} \chi_0  ~&=&~ \{m_{1}, m_{2}; ~ -\ha (2m_{2}+m_{1})\} \\
\chi_2 ~&=&~ \{ 3m_{2}+m_{1}, -m_{2}  ; ~ -\ha (m_{2}+m_{1}) \}, \quad \L_{2} = \L_0-m_{2}\a_{2}  \nn\\
\chi_1 ~&=&~ \{ -m_{1}, m_{2}+m_{1} ;~  -\ha (2m_{2}+m_{1}) \}, \quad \L_{1} = \L_0-m_{1}\a_{1}  \nn\\
\chi_{21} ~&=&~ \{3m_{2}+2m_{1}, -m_{2}-m_{1}   ;~ -\ha m_{2}  \},  \quad \L_{21} = \L_{2}- m_{1} \a_{3}\nn\\
\chi_{12} ~&=&~ \{ -3m_{2}-m_{1} , 2m_{2}+m_{1} ; ~ -\ha (m_{2}+m_{1})  \}, \quad \L_{12} = \L_{1}- m_{2} \a_{5} \nn\\
\chi_{212} ~&=&~ \{ 3m_{2}+2m_{1}, -2m_{2}-m_{1}   ;~ \ha m_{2}   \},  \quad \L_{212} = \L_0- (m_1+3m_2)\a_{3} \nn\\ % (\L_{21}- m_{2} \a_{6}\nn\\
\chi_{121} ~&=&~ \{ -3m_{2}-2m_{1}, 2m_{2}+m_{1}  ;~ -\ha m_{2}  \}, \quad \L_{121} = \L_0- (m_1+m_2)\a_{5} \nn\\ %\L_{12}- m_{1} \a_{4}\nn\\
\chi_{2121} ~&=&~ \{3m_{2}+m_{1}, -2m_{2}-m_{1}    ;~ \ha (m_{2}+m_{1})   \},  \quad \L_{2121} = \L_{212}- m_{1} \a_{4}\nn\\
\chi_{1212} ~&=&~ \{ -3m_{2}-2m_{1}, m_{2}+m_{1}  ;~ \ha m_{2}  \}, \quad \L_{1212} = \L_{121}- m_{2} \a_{6}\nn\\
\chi_{21212} ~&=&~ \{m_{1}, -m_{2}-m_{1}   ;~  \ha (2m_{2}+m_{1})  \},  \quad \L_{21212} = \L_0- (m_1+2m_2)\a_{6} \nn\\  %\L_{2121}- m_{2} \a_{5}\nn\\
\chi_{12121} ~&=&~ \{ -3m_{2}-m_{1}, m_{2}  ; ~ \ha (m_{2}+m_{1}) \}, \quad \L_{12121} = \L_0- (2m_1+3m_2)\a_{4} \nn\\   %\L_{1212}- m_{1} \a_{3}\nn\\
\chi_{121212} ~&=&~ \{-m_{1}, -m_{2}   ; ~ \ha (2m_{2}+m_{1})  \} =\chi_{212121} = (s_1s_2)^3\cdot \L_0 = (s_2s_1)^3\cdot \L_0 , \nn\ee
%&&  \L_{121212} = \L_{12121}- m_{1} \a_{1}= \L_{21212}- m_{2} \a_{2}\nn   \\
%&& \L_{21} = \L_{2}- (3m_{2}+m_{1}) \a_{1} \nn\\
%&&   \L_{212} = \L_{12} - (3m_{2}+2m_{1}) \a_{1} \nn\\
%&&   \L_{2121} = \L_{121} - (3m_{2}+2m_{1}) \a_{1} \nn\\
%&&   \L_{21212} = \L_{1212} - (3m_{2}+m_{1}) \a_{1} \nn\\
%&&\L_{12} = \L_{1}- (m_{2}+m_{1}) \a_{2} \nn\\
%&&\L_{121} = \L_{21}- (2m_{2}+m_{1}) \a_{2} \nn\\
%&&\L_{1212} = \L_{212}- (2m_{2}+m_{1}) \a_{2} \nn\\
%&&\L_{12121} = \L_{2121}- (m_{2}+m_{1}) \a_{2} \nn
where we have included as third entry also the parameter ~$c ~=~ -\ha (2m_{2}+m_{1})$, related to the Harish-Chandra parameter of the highest root
(recalling that ~$m_{\a_{6}} =2m_{2}+m_{1} $). It is also related to the conformal weight ~$d = \trh +c$.

Using  this labelling  the signatures may be
given in the following pair-wise manner:
\eqnn{enumpm} \chi^\pm_0  ~&=&~ \{\mp m_{1},\mp m_{2}; ~ \pm\ha (2m_{2}+m_{1})\}  \\
\chi^\pm_{2} ~&=&~ \{ \mp(3m_{2}+m_{1}) ,  \pm m_{2} ; ~ \pm\ha (m_{2}+m_{1}) \},   \nn\\
\chi^\pm_{1} ~&=&~ \{ \pm m_{1} , \mp (m_{2}+m_{1}) ;~  \pm\ha (2m_{2}+m_{1}) \},   \nn\\
\chi^\pm_{12} ~&=&~ \{ \mp (3m_{2}+2m_{1}) , \pm (m_{2}+m_{1})  ;~ \pm \ha m_{2}  \}  \nn\\
\chi^\pm_{21} ~&=&~ \{ \pm (3m_{2}+m_{1}) ,\mp (2m_{2}+m_{1}) ; ~ \pm\ha (m_{2}+m_{1})  \} \nn\\
\chi^\pm_{121} ~&=&~ \{  \mp (3m_{2}+2m_{1}) ,\pm (2m_{2}+m_{1})  ;~ \mp\ha m_{2}   \} ,   \nn
 \ee
where ~$\chi^-_{...} =  \chi_{...}$~ from \eqref{enum}, ~$\chi^+_0 =\chi_{121212}$,
~$\chi^+_1~=~ \chi_{21212}$, ~ $\chi^+_2~=~  \chi_{12121}$, ~$\chi^+_{12}~=~ \chi_{2121}$,
 ~$\chi^+_{21} ~=~ \chi_{1212}$,  ~$\chi^+_{121} ~=~\chi_{212}$~.

%##

The ERs in the multiplet are related also by intertwining integral
  operators  introduced in \cite{KnSt}. These operators are defined
for any ER,   the general action in our situation
being: \eqnn{knast}  && G_{KS} ~:~ \cc_\chi ~ \llr ~ \cc_{\chi'} \
,\nn\\ &&\chi ~=~ [\,    n_{1},n_{2}
\,;\, c\, ] \ , \qquad \chi' ~=~ [\,
    -n_{1},-n_{2}  \,;\, -c\, ] . \ee
This action is consistent with the parametrization in \eqref{enumpm}.

 The main multiplets are given explicitly in Fig.~1.
The pairs ~$\chi^\pm$~ are symmetric w.r.t.  the bullet in the middle of the picture - this symbolizes the Weyl symmetry realized by
the Knapp-Stein operators \eqref{knast}:\\ $G^\pm ~:~ \cc_{\chi^\mp}
\lra \cc_{\chi^\pm}\,$.

 \bigskip

Some comments are in order.

Matters are arranged so that in every multiplet only the ER with
signature ~$\chi_0^-$~ contains a finite-dimensional nonunitary
subrepresentation in  a finite-dimensional subspace ~$\ce$. The
latter corresponds to the finite-dimensional   irrep of ~$G_{2 (2)}$~ with
signature ~$[  m_{1},m_{2} ]$.   The subspace ~$\ce$~ is annihilated by the
operators ~$G^+\,$,\ $\cd^{m_1}_{\a_1}$,\ $\cd^{m_2}_{\a_2}$\ and is the image of the operator ~$G^-\,$.\\
 When both ~$m_i=1$~ then ~$\dim\,\ce = 1$, and in that case
~$\ce$~ is also the trivial one-dimensional UIR of the whole algebra
~$\cg$. Furthermore in that case the conformal weight is zero:
~$d=\trh+c=\trh-\ha (2m_{2}+m_{1})_{\vert_{m_i=1}}=0$.

In the conjugate ER ~$\chi_0^+$~ there is a unitary discrete series
 representation (according to the Harish-Chandra criterion \cite{Har}) % that the Harish-Chandra parameters are positive for the
%$\ck$-compact roots $\a_3,\a_5$)
in  an infinite-dimen\-sional subspace  $\tcd_0$ with conformal weight
~$d=\trh+c=\trh+\ha (2m_{2}+m_{1}) ~=~ 3, \sevha,4,...$. It
is annihilated by the operator  $G^- $,\ and is in the intersection of the images of the
operators  $G^+ $ (acting from $\chi_0^-$), $\cd^{m_1}_{\a_1}$ (acting from $\chi_1^+$), $\cd^{m_2}_{\a_2}$ (acting from $\chi_2^+$).

\bigskip

\nt {\bf Remark:}~ In paper \cite{Dobg2} were considered also the following multiplets for $\cgc$ which are not interesting for the real form $\cg$.
 Fix $k=1,\ldots,6$.
  Then there are Verma modules multiplets %of   subtype $N_k$ are
  parametrized by the natural number $m_k$, so that $m_j \notin \bbn$ for $j\neq k$,
and  given as follows:
$ %\eqn{typa}
 V^{\L_k} ~\longrightarrow~   V^{\L_k - m_k\a_k} \  . $~\dia %\eeq %&\typa a\cr

\subsection{Reduced multiplets}

There are  two   reduced multiplets $M_k$, $k=1,2$, which may be obtained by setting the parameter $m_k=0$.

The reduced multiplet $M_{1}$ contains  six   GVMs (ERs).
 Their signatures are given as follows:
  \eqnn{enumpm1} \chi'^\pm_0  ~&=&~ \{  0,\mp m_{2}; ~ \pm m_{2}\} =\chi'^\pm_{1} \\
 \chi'^\pm_{2} ~&=&~ \{ \mp 3m_{2} ,  \pm m_{2} ; ~ \pm\ha m_{2} \}=\chi'^\pm_{12} ,   \nn\\
%\chi'^\p0 ~&=&~ \{ \pm 0 , \mp (m_{2}+0) ;~  \pm\ha (2m_{2} \},   \nn\\
%\chi'^\pm_{12} ~&=&~ \{ \mp 3m_{2} , \pm m_{2}  ;~ \pm \ha m_{2}  \}  \nn\\
\chi'^\pm_{21} ~&=&~ \{ \pm 3m_{2} ,\mp 2m_{2}  ; ~ \pm\ha m_{2}  \} =\chi'^\mp_{121} \nn
%\chi'^\pm_{121} ~&=&~ \{  \mp 3m_{2} ,\pm 2m_{2}  ;~ \mp\ha m_{2}   \} ,   \nn
 \ee

 The \idos\ of the multiplet are given explicitly as follows, cf. \eqref{lido}:
 \eqn{enumr11}  \L'^{-}_0 ~{\cd^{m_{2}}_{\a_2}  \atop \longrightarrow } ~ \L'^-_2 ~{ \left( \cd^{m_{2}}_{\a_1}\right)^3  \atop \longrightarrow } ~
\L'^-_{21} ~ {\left( \cd^{m_{2}}_{\a_2}\right)^2  \atop \longrightarrow } ~ \L'^+_{21} { \left( \cd^{m_{2}}_{\a_1}\right)^3  \atop \longrightarrow }
\L'^+_{2} ~{\cd^{m_{2}}_{\a_2}  \atop \longrightarrow } ~
 ~\L'^+_0 \eeq
Note a peculiarity on the map from $\L'^-_{21}$ to $\L'^+_{21}$ - it is a degeneration of the
corresponding $G^+$ KS operator. In addition it is a part of the chain degeneration
of the $G^+$ KS operators   from $\L'^-_{2}$ to $\L'^+_{2}$
 and   from $\L'^-_{0}$ to $\L'^+_{0}$.
 Thus, the diagram  may be represented also as:
\eqn{enumrq1a} \begin{matrix} \L'^-_0 & ~{\cd^{m_{2}}_{\a_2}  \atop \longrightarrow }  & \L'^-_2 & ~{ \left( \cd^{m_{2}}_{\a_1}\right)^3  \atop \longrightarrow } ~
 & \L'^-_{21} \cr
&&& \cr
\updownarrow & &\updownarrow && \uparrow G^- ~ \downarrow  \left( \cd^{m_{2}}_{\a_2}\right)^2   \cr
&&&\cr
\L'^+_0 & { \cd^{m_{2}}_{\a_2}    \atop \longleftarrow } & \L'^+_2 & {  \left( \cd^{m_{2}}_{\a_1}\right)^3  \atop \longleftarrow } & \L'^+_{21}
\end{matrix}
 \eeq
The ER ~$\chi'^+_0$~ contains  a unitary discrete series
 representation  in  an infinite-dimen\-sional subspace  $\tcd_1$ with conformal weight
~$d=\trh+c= \fha,3,\sevha,...$. % (here $c_d=1$).
It is in the intersection of the images of the
operators  $G^+ $ (acting from $\chi'^-_0$) and  $\cd^{m_2}_{\a_2}$ (acting from $\chi'^+_2$).
It is different from the case $\tcd_0$ even when the conformal weights coincide.

Note also that the discrete series representation in $\tcd_1$ may be
obtained as a subrepresentation when inducing from maximal parabolic
$\cp_1$, see corresponding section below.

\bigskip

The   reduced multiplet $M_{2}$ contains  six   GVMs (ERs).
 Their signatures are given as follows:

   \eqnn{enumpm2} \chi''^\pm_0  ~&=&~ \{\mp m_{1}, 0; ~ \pm\ha m_{1}\} =\chi''^\pm_{2}  \\
%\chi''^\pm_{2} ~&=&~ \{ \mp m_{1} ,   0 ; ~ \pm\ha m_{1} \},   \nn\\
\chi''^\pm_{1} ~&=&~ \{ \pm m_{1} , \mp m_{1} ;~  \pm\ha m_{1} \} =\chi''^\pm_{21},   \nn\\
\chi''^\pm_{12} ~&=&~ \{ \mp 2m_{1} , \pm m_{1}  ;~  0  \}  =\chi''^\pm_{121}\nn
%\chi''^\pm_{21} ~&=&~ \{ \pm m_{1} ,\mp m_{1} ; ~ \pm\ha m_{1}  \} \nn\\
%\chi''^\pm_{121} ~&=&~ \{  \mp 2m_{1} ,\pm m_{1}  ;~  0   \} ,   \nn
 \ee

 The \idos\ of the  multiplet  are given explicitly as follows:
 \eqn{enumr21}  \L''^{-}_0 ~{\cd^{m_{1}}_{\a_1} \atop \longrightarrow } ~ \L''^-_1 ~{  \cd^{m_{1}}_{\a_2}  \atop \longrightarrow } ~
\L''^-_{12} ~{\left(\cd^{m_{1}}_{\a_1}\right)^2 \atop \longrightarrow } ~ \L''^+_{12} ~{ \cd^{m_{1}}_{\a_2}   \atop \longrightarrow } ~
\L''^+_1 ~{ \cd^{m_{1}}_{\a_1}   \atop \longrightarrow } ~\L''^+_0 \eeq

Also here we note a peculiarity similar to the previous case. The map from $\L''^-_{12}$ to $\L''^+_{12}$  is    degeneration of the
corresponding $G^+$ KS operator. In addition it is a part of the chain degeneration
of the $G^+$ KS operators from $\L''^-_{1}$ to $\L''^+_{1}$ and from $\L''^-_{0}$ to $\L''^+_{0}$. Thus, the diagram  may be represented also as:
\eqn{enumr21a} \begin{matrix} \L''^-_0 & {\cd^{m_{1}}_{\a_1} \atop \longrightarrow } & \L''^-_1 &{ \cd^{m_{1}}_{\a_2}   \atop \longrightarrow } &
\L''^-_{12} \cr
&&& \cr
\updownarrow & &\updownarrow && \uparrow G^- ~\downarrow \left(\cd^{m_{1}}_{\a_1}\right)^2  \cr
&&&\cr
\L''^+_0 & { \cd^{m_{1}}_{\a_1}   \atop \longleftarrow } & \L''^+_1 & { \cd^{m_{1}}_{\a_2}   \atop \longleftarrow } & \L''^+_{12}
\end{matrix}
 \eeq

The ER ~$\chi''^+_0$~ contains  a unitary discrete series
 representation  in  an infinite-dimen\-sional subspace  $\tcd_2$ with conformal weight
~$d=\trh+c= 2,\fha,3,\sevha,...$. % (here $c_d=\ha$).
It
is in the intersection of the images of the
operators  $G^+ $ (acting from $\chi''^-_0$) and  $\cd^{m_1}_{\a_1}$ (acting from $\chi''^+_1$).
It is different from the cases $\tcd_0,\tcd_1$ even when the conformal weights coincide.

Note also that the discrete series representation in $\tcd_2$ may be
obtained as a subrepresentation when inducing from maximal parabolic
$\cp_2$, see corresponding section below.

 The  multiplets in this section are shared with $G_2$ as Verma modules multiplets and were given in \cite{Dobg2} but without the weights of the singular vectors, the KS operators
 and the identification of the discrete series.

\bigskip

\section{Induction from maximal parabolics}

 As stated in Section 1 in order to obtain all possible \idos\ we should consider induction from all parabolics, yet taking into account the possible overlaps.

\subsection{Main multiplets when inducing from  $\cp_1$}

When inducing from the maximal parabolic $\cp_1= \cm_{1} \oplus \ca_1 \oplus \cn_1$ there is one $\cm_{1}$-compact root, namely, $\a_1$.
We take again the Verma module with $\L_{HC}=\L^{1-}_0$. We take ~$\chi^{1-}_0 = \chi_{HC}$.
The GVM  $\L^{1-}_0$  has one embedded GVM  with HW ~$\L^{1-}_2 = \L^{1-}_0-m_{2}\a_{2}$, ~$m_{2}\in\bbn$.
Altogether, the main multiplet in this case includes the same number of  ERs/GVMs as in \eqref{enum}, so we use the same notation only adding super index 1, namely
 \eqnn{enumpm1a} \chi^{1\pm}_0  ~&=&~ \{\mp m_{1},\mp m_{2}; ~ \pm\ha (2m_{2}+m_{1})\}  \\
\chi^{1\pm}_{2} ~&=&~ \{ \mp(3m_{2}+m_{1}) ,  \pm m_{2} ; ~ \pm\ha (m_{2}+m_{1}) \},   \nn\\
\chi^{1\pm}_{1} ~&=&~ \{ \pm m_{1} , \mp (m_{2}+m_{1}) ;~  \pm\ha (2m_{2}+m_{1}) \},   \nn\\
\chi^{1\pm}_{12} ~&=&~ \{ \mp (3m_{2}+2m_{1}) , \pm (m_{2}+m_{1})  ;~ \pm \ha m_{2}  \}  \nn\\
\chi^{1\pm}_{21} ~&=&~ \{ \pm (3m_{2}+m_{1}) ,\mp (2m_{2}+m_{1}) ; ~ \pm\ha (m_{2}+m_{1})  \} \nn\\
\chi^{1\pm}_{121} ~&=&~ \{  \mp (3m_{2}+2m_{1}) ,\pm (2m_{2}+m_{1})  ;~ \mp\ha m_{2}   \} ,   \nn
 \ee
  In addition, in order to avoid coincidence with \eqref{enumpm}  we must impose in \eqref{enumpm1a} the
conditions: ~$m_{1}\notin\bbn$, ~$m_{1}\notin\bbn/2$.
%The ER ~$\chi^{1+}_{1}$~ contain discrete series representation  since all $\cm_{2}$-noncompact HC parameters are negative.

What is peculiar is that the ERs/GVMs of the main miltiplet \eqref{enumpm1a} actually consists  of three submultiplets with intertwining diagrams as follows:
   \eqna{diagp1}  && \begin{matrix}
  \L^{1-}_0 &{\cd^{m_{2}}_{\a_{2}}  \atop\llr }&\L^{1-}_2  \cr
  &&\cr
  \updownarrow && \updownarrow \cr
  &&\cr
  \L^{1+}_0 &{\cd^{m_{2}}_{\a_{2}}  \atop\lla }&\L^{1+}_2  \cr
  \end{matrix} \quad {\rm subtype~ (A_1)}\\
  && \phantom{ \L^{1+}_1} \nn\\ %\\[10pt]%  \vspace{1cm}
  &&\begin{matrix}
  \L^{1-}_1 & {\cd^{m_{2}}_{\a_{5}}  \atop\llr }&\L^{1-}_{21}  \cr
  &&\cr
  \updownarrow && \updownarrow \cr
  &&\cr
  \L^{1+}_1 &{\cd^{m_{2}}_{\a_{5}}  \atop\lla }&\L^{1+}_{21}  \cr
  \end{matrix} \quad {\rm subtype~ (B_1)}
  \\
  && \phantom{ \L^{1+}_1} \nn\\
  &&\begin{matrix}
  \L^{1-}_{12} & {\cd^{m_{2}}_{\a_{6}}  \atop\llr }&\L^{1+}_{121}  \cr
  &&\cr
  \updownarrow && \updownarrow \cr
  &&\cr
  \L^{1+}_{12} &{\cd^{m_{2}}_{\a_{6}}  \atop\lla }&\L^{1-}_{121}  \cr
  \end{matrix} \quad {\rm subtype~ (C_1)}
  \eena

  Next we   relax in \eqref{enumpm1a} one of the conditions, namely, we allow ~$m_{1}\in \bbn/2$,  still keeping
  ~$m_{2}\in\bbn$,   ~$m_{1}\notin\bbn$. This changes the diagram of subtype ($C_1$), \rf{diagp1}c, as given in Fig. 2a.
    %\fig{}{diag-g2-pa.eps}{20cm}

  \subsection{Main multiplets when induction from  $\cp_2$}

This case is  partly dual to the previous one.
When inducing from the maximal parabolic $\cp_2= \cm_{2} \oplus \ca_2 \oplus \cn_2$ there is one $\cm_{2}$-compact root, namely, $\a_2$.
We take again the Verma module with $\L_{HC}=\L^{2-}_0$. We take ~$\chi^{2-}_0 = \chi_{HC}$.
The GVM  $\L^{2-}_0$  has one embedded GVM  with HW ~$\L^{2-}_1 = \L^{2-}_0-m_{1}\a_{1}$, ~$m_{1}\in\bbn$.
Altogether, the main multiplet in this case includes the same number of  ERs/GVMs as in \eqref{enum}, so we use the same notation only adding super index 2, namely
 \eqnn{enumpm2a} \chi^{2\pm}_0  ~&=&~ \{\mp m_{1},\mp m_{2}; ~ \pm\ha (2m_{2}+m_{1})\}  \\
\chi^{2\pm}_{2} ~&=&~ \{ \mp(3m_{2}+m_{1}) ,  \pm m_{2} ; ~ \pm\ha (m_{2}+m_{1}) \},   \nn\\
\chi^{2\pm}_{1} ~&=&~ \{ \pm m_{1} , \mp (m_{2}+m_{1}) ;~  \pm\ha (2m_{2}+m_{1}) \},   \nn\\
\chi^{2\pm}_{12} ~&=&~ \{ \mp (3m_{2}+2m_{1}) , \pm (m_{2}+m_{1})  ;~ \pm \ha m_{2}  \}  \nn\\
\chi^{2\pm}_{21} ~&=&~ \{ \pm (3m_{2}+m_{1}) ,\mp (2m_{2}+m_{1}) ; ~ \pm\ha (m_{2}+m_{1})  \} \nn\\
\chi^{2\pm}_{121} ~&=&~ \{  \mp (3m_{2}+2m_{1}) ,\pm (2m_{2}+m_{1})  ;~ \mp\ha m_{2}   \} ,   \nn
 \ee

In addition, in order to avoid coincidence with \eqref{enumpm}  we must impose in \eqref{enumpm2a} the
conditions: ~$m_{2}\notin\bbn$, ~$m_{2}\notin\bbn/2$, ~$m_{2}\notin\bbn/3$.

Similarly to the $\cp_1$ case  the ERs/GVMs of the main miltiplet \eqref{enumpm2a} actually consists  of three submultiplets with intertwining diagrams as follows:
   \eqna{diagp2}  && \begin{matrix}
  \L^{2-}_0 &{\cd^{m_{1}}_{\a_{1}}  \atop\llr }&\L^{2-}_1  \cr
  &&\cr
  \updownarrow && \updownarrow \cr
  &&\cr
  \L^{2+}_0 &{\cd^{m_{1}}_{\a_{1}}  \atop\lla }&\L^{2+}_1  \cr
  \end{matrix} \quad {\rm subtype~ (A_2)}\\
  && \phantom{ \L^{2+}_1} \nn\\ %\\[10pt]%  \vspace{1cm}
  &&\begin{matrix}
  \L^{2-}_2 & {\cd^{m_{1}}_{\a_{3}}  \atop\llr }&\L^{2-}_{12}  \cr
  &&\cr
  \updownarrow && \updownarrow \cr
  &&\cr
  \L^{2+}_2 &{\cd^{m_{1}}_{\a_{3}}  \atop\lla }&\L^{2+}_{12}  \cr
  \end{matrix} \quad {\rm subtype~ (B_2)}
  \\
  && \phantom{ \L^{2+}_1} \nn\\
  &&\begin{matrix}
  \L^{2-}_{21} & {\cd^{m_{1}}_{\a_{4}}  \atop\llr }&\L^{2-}_{121}  \cr
  &&\cr
  \updownarrow && \updownarrow \cr
  &&\cr
  \L^{2+}_{21} &{\cd^{m_{1}}_{\a_{4}}  \atop\lla }&\L^{2+}_{121}  \cr
  \end{matrix} \quad {\rm subtype~ (C_2)}
  \eena

\bigskip

  Next we   relax in \eqref{enumpm2a} one of the conditions, namely, we allow ~$m_{2}\in \bbn/2$,  still keeping
  ~$m_{2}\notin\bbn$, ~$m_{2}\notin\bbn/3$. This changes the diagram of subtype ($C_2$), \rf{diagp2}c, as given in Fig. 2b.\\
 In this case the ER ~$\chi_0^{2+}$~ contains  a %unitary discrete series
subrepresentation  in  an infinite-dimen\-sional subspace  $\tcd'_0$ with conformal weight
~$d=\trh+c=  \trh + m_2+ \ha m_{1} = \fha,3,\sevha,...$. It
is in the intersection of the images of the
operators  $G^+ $ (acting from $\chi_0^{2-}$) and  $\cd^{m_1}_{\a_1}$ (acting from $\chi_1^{2+}$).

  Next we   relax in \eqref{enumpm2a}  another condition, namely, we allow ~$m_{2}\in \bbn/3$,  still keeping
  ~$m_{2}\notin\bbn$,   ~$m_{2}\notin\bbn/2$. This changes the diagrams of subtypes ($B_1$) and ($C_1$) combining them as  given in Fig. 2c.
   % \fig{}{diag-g2-pb.eps}{20cm}

  \bigskip

 \subsection{Reduced multiplets}

There are  two reduced multiplets $M_{1k}$, $k=1,2$, which may be obtained by setting the parameter $m_k=0$ when inducing from $\cp_1$,
and two reduced multiplets $M_{2k}$, $k=1,2$, which may be obtained by setting the parameter $m_k=0$ when inducing from $\cp_2$.
 %What is peculiar that   these two cases are interrelated.

In case $M_{11}$ ($m_{1}=0$, $m_{2}\in\bbn$) the reduced multiplet has  six GVMs:
\eqnn{enumpm11} \chi^{1\pm}_0  ~&=&~ \{ 0,\mp m_{2}; ~ \pm m_{2}\}  =\chi^{1\pm}_{1} \\
\chi^{1\pm}_{2} ~&=&~ \{ \mp 3m_{2} ,  \pm m_{2} ; ~ \pm\ha m_{2} \}=\chi^{1\pm}_{12}  \nn\\
%\chi^{1\pm}_{1} ~&=&~ \{  0 , \mp m_{2} ;~  \pm m_{2} \},   \nn\\
%\chi^{1\pm}_{12} ~&=&~ \{ \mp 3m_{2} , \pm m_{2}  ;~ \pm \ha m_{2}  \}  \nn\\
\chi^{1\pm}_{21} ~&=&~ \{ \pm 3m_{2} ,\mp 2m_{2} ; ~ \pm\ha m_{2}  \} =\chi^{1\pm}_{121}\nn
%\chi^{1\pm}_{121} ~&=&~ \{  \mp 3m_{2} ,\pm 2m_{2}  ;~ \mp\ha m_{2}   \} ,   \nn
 \ee
Note that thus reduced multiplet coincides with the reduced multiplet $M_{1}$ when inducing from the minimal parabolic, cf.
\eqref{enumpm1}.
 The  \idos\ are correspondingly given by a recombination of the three   submultiplets  from \eqref{diagp1} and the resulting
 diagram coincides with the one in \eqref{enumrq1a}.

\bigskip

In case $M_{12}$ (when $m_{2}=0$) the reduced multiplet has  six GVMs:
  \eqnn{enumpm12} \chi^{1\pm}_0  ~&=&~ \{\mp m_{1}, 0; ~ \pm\ha m_{1}\} =\chi^{1\pm}_{2} \\
%\chi^{1\pm}_{2} ~&=&~ \{ \mp m_{1} ,   0 ; ~ \pm\ha m_{1} \},   \nn\\
\chi^{1\pm}_{1} ~&=&~ \{ \pm m_{1} , \mp m_{1} ;~  \pm\ha m_{1} \}=\chi^{1\pm}_{21}   \nn\\
\chi^{1\pm}_{12} ~&=&~ \{ \mp 2m_{1} , \pm m_{1}  ;~ 0  \} =\chi^{1\mp}_{121} \nn
%\chi^{1\pm}_{21} ~&=&~ \{ \pm m_{1} ,\mp m_{1} ; ~ \pm\ha m_{1}  \} \nn\\
%\chi^{1\pm}_{121} ~&=&~ \{  \mp 2m_{1} ,\pm m_{1}  ;~ 0   \} ,   \nn
 \ee

As for the main multiplet we first consider the subcase
~$m_{1}\notin\bbn$, ~$m_{1}\notin\bbn/2$. Again as in \eqref{diagp1} we have three submultiplets, however the submultiplets
$(A_1)$, $(B_1)$, $(C_1)$, are replaced by
  the KS related doublets $\chi^{1\pm}_0$, $\chi^{1\pm}_1$, $\chi^{1\pm}_{12}$.

Next we consider the subcase ~$m_{1}\notin\bbn$, ~$m_{1}\in\bbn/2$. As in the first subcase we have the three KS related doublets,
yet for the doublet $\chi^{1\pm}_{12}$ the $G^+$ KS operator degenerates to the \ido\ ~$\cd^{2m_{1}}_{\a_{1}}$, (compare with Fig. 2a).

\bigskip

In case $M_{21}$ ($m_{1}=0$) the reduced multiplet has  six GVMs with signatures:
\eqnn{enumpm21} \chi^{2\pm}_0  ~&=&~ \{ 0,\mp m_{2}; ~ \pm m_{2}\}=\chi^{2\pm}_{1}  \\
\chi^{2\pm}_{2} ~&=&~ \{ \mp 3m_{2} ,  \pm m_{2} ; ~ \pm\ha m_{2} \}=\chi^{2\pm}_{12}   \nn\\
%\chi^{2\pm}_{1} ~&=&~ \{ 0, \mp m_{2} ;~  \pm m_{2} \},   \nn\\
%\chi^{2\pm}_{12} ~&=&~ \{ \mp 3m_{2} , \pm m_{2}  ;~ \pm \ha m_{2}  \}  \nn\\
\chi^{2\pm}_{21} ~&=&~ \{ \pm 3m_{2} ,\mp 2m_{2} ; ~ \pm\ha m_{2}  \} =\chi^{2\mp}_{121}\nn
%\chi^{2\pm}_{121} ~&=&~ \{  \mp 3m_{2} ,\pm 2m_{2}  ;~ \mp\ha m_{2}   \} ,   \nn
 \ee

As for the main multiplet we first consider the subcase
~$m_{2}\notin\bbn$, ~$m_{2}\notin\bbn/2$, ~$m_{2}\notin\bbn/3$.
Again as in \eqref{diagp1} we have three submultiplets, however the submultiplets
$(A_2)$, $(B_2)$, $(C_2)$, are replaced by
  the KS related doublets $\chi^{2\pm}_0$, $\chi^{2\pm}_2$, $\chi^{2\pm}_{21}$.

Next we consider the subcase  ~$m_{2}\notin\bbn$, ~$m_{2}\in\bbn/2$, ~$m_{2}\notin\bbn/3$.
As in the first subcase we have the three KS related doublets,
yet for the doublet $\chi^{2\pm}_{21}$ the $G^+$ KS operator degenerates to the \ido\ ~$\cd^{2m_{2}}_{\a_{2}}$, compare with Fig. 2c.\\
%Conjecture:
The ER ~$\chi_0^{2+}$~ contains  a %unitary discrete series
subrepresentation  in  an infinite-dimen\-sional subspace  $\tcd'_1$ with conformal weight
~$d=\trh+c=  \trh + m_2  =2,3,4,...$.% ; (here $c_d = \ha$).
It is the   image of the KS
operator  $G^+ $ (acting from $\chi_0^{2-}$).

Next we consider the subcase  ~$m_{2}\notin\bbn$, ~$m_{2}\notin\bbn/2$, ~$m_{2}\in\bbn/3$.
As for the main multiplet the submultiplets
 corresponding to $(B_2)$, $(C_2)$ are combined.
 Here the result is a quartet (compare with Fig. 2c):
 \eqn{diagp12} \begin{matrix} \chi^{2-}_2 & {\cd^{3m_2}_{\a_{1}}\atop \llr} & \chi^{2-}_{21} \cr
 &&\cr
 \updownarrow &&  \updownarrow \cr
 &&\cr
\chi^{2+}_2 & {\cd^{3m_2}_{\a_{1}}\atop \lla} & \chi^{2+}_{21}
\end{matrix}
\eeq

\bigskip

In case $M_{22}$ ($m_{2}=0$, $m_{1}\in\bbn$) the reduced multiplet has  six GVMs with signatures:
 \eqnn{enumpm22} \chi^{2\pm}_0  ~&=&~ \{\mp m_{1}, 0; ~ \pm\ha m_{1}\} =\chi^{2\pm}_{2} \\
%\chi^{2\pm}_{2} ~&=&~ \{ \mp m_{1} ,   0 ; ~ \pm\ha m_{1} \},   \nn\\
\chi^{2\pm}_{1} ~&=&~ \{ \pm m_{1} , \mp m_{1} ;~  \pm\ha m_{1} \}=\chi^{2\pm}_{21}   \nn\\
\chi^{2\pm}_{12} ~&=&~ \{ \mp 2m_{1} , \pm m_{1}  ;~  0  \} =\chi^{2\pm}_{121} \nn
%\chi^{2\pm}_{21} ~&=&~ \{ \pm m_{1} ,\mp m_{1}  ; ~ \pm\ha m_{1}  \} \nn\\
%\chi^{2\pm}_{121} ~&=&~ \{  \mp 2m_{1} ,\pm m_{1}  ;~  0   \} ,   \nn
 \ee
Note that thus reduced multiplet coincides with the reduced multiplet $M_{2}$ when inducing from the minimal parabolic, cf.
\eqref{enumpm2}.
 The  \idos\ are correspondingly given by a recombination of the three   submultiplets  from \eqref{diagp2} and the resulting diagram coincides with the one in \eqref{enumr21}.

\np

\section*{Acknowledgments.}

\nt The author acknowledges partial support from Bulgarian NSF Grant
DN-18/1.

\newpage

%$$\includegraphics[width=20cm] {diag-sost8}%$$
%%$$\includegraphics[bb= 0 -10 500 120, width=10cm] {diag-sost8}%$$

 \fig{}{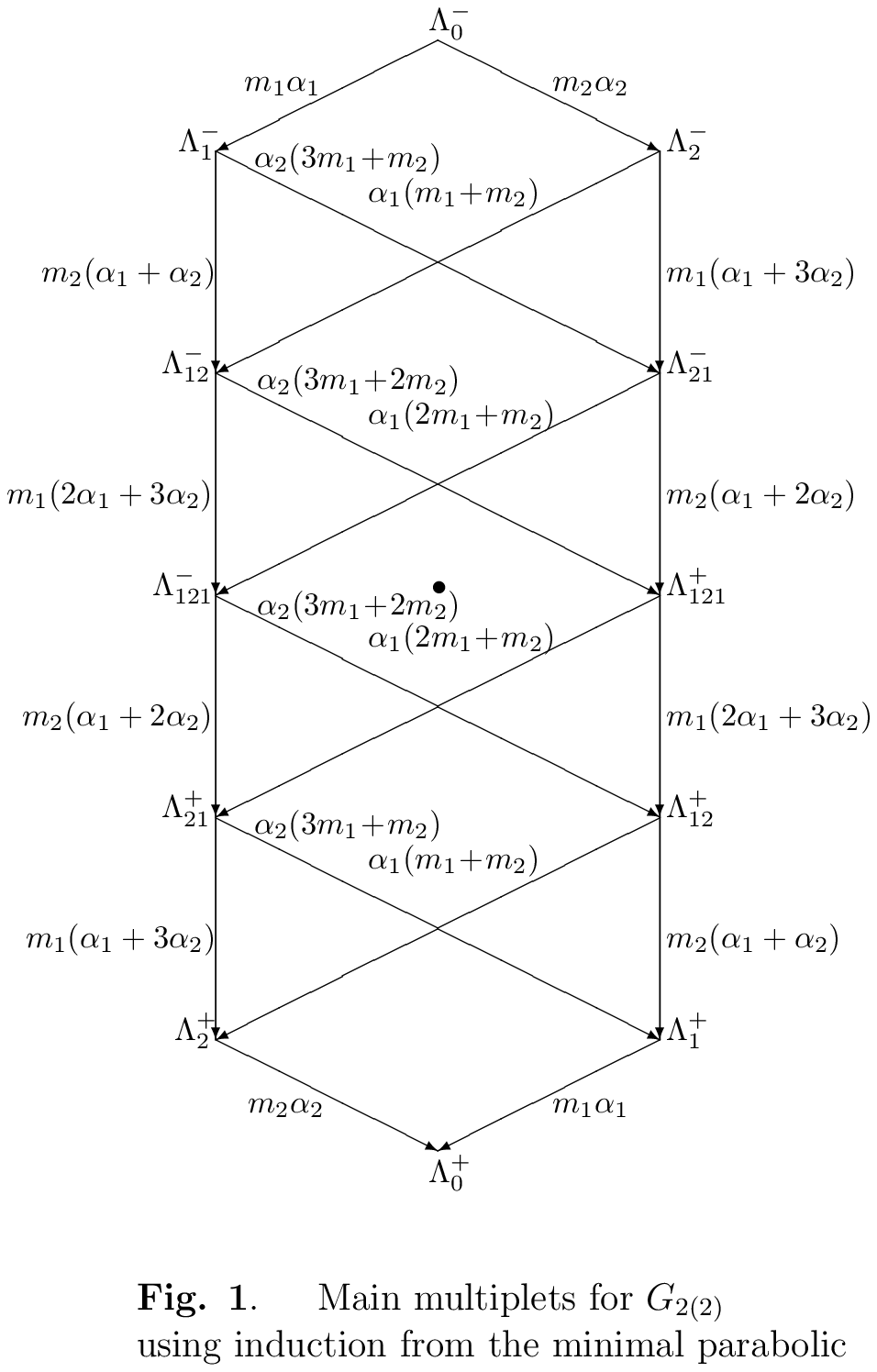}{20cm}

 \fig{}{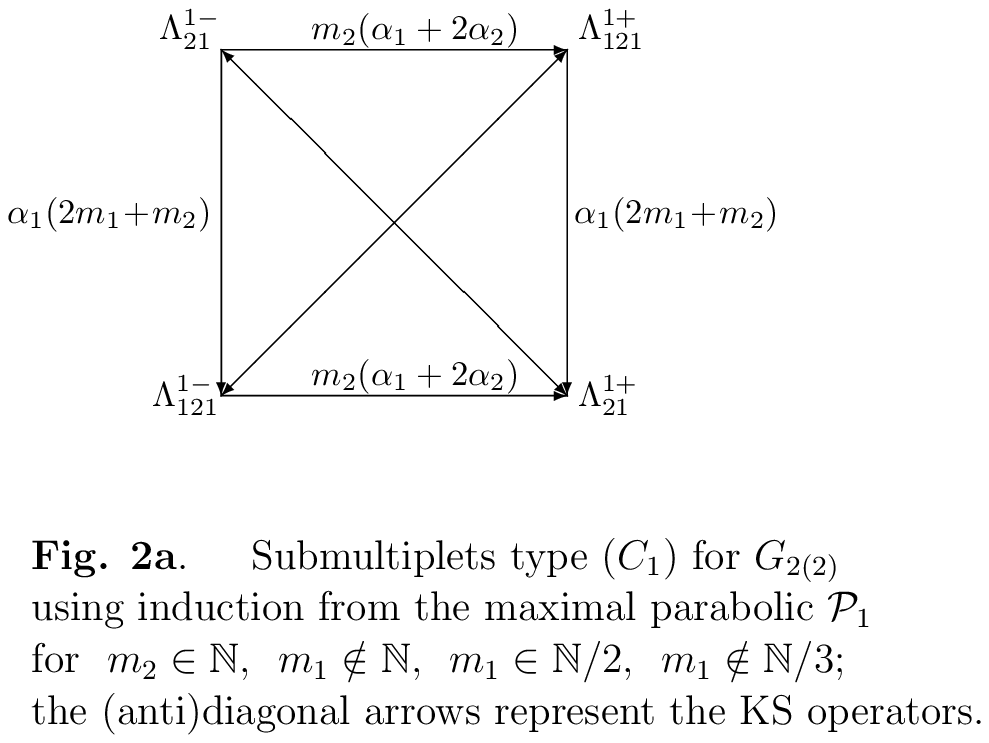}{20cm}

 \fig{}{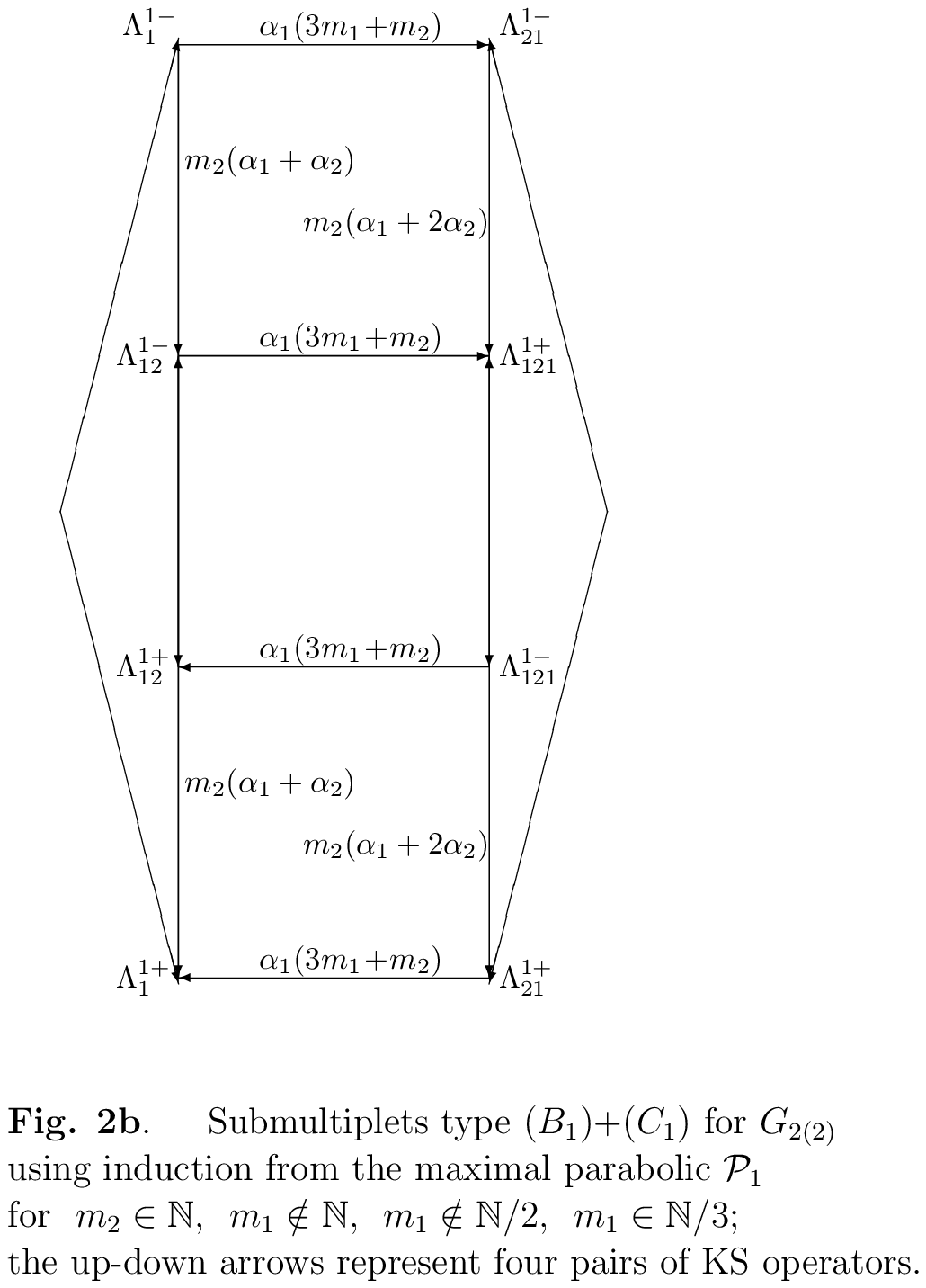}{20cm}

 \fig{}{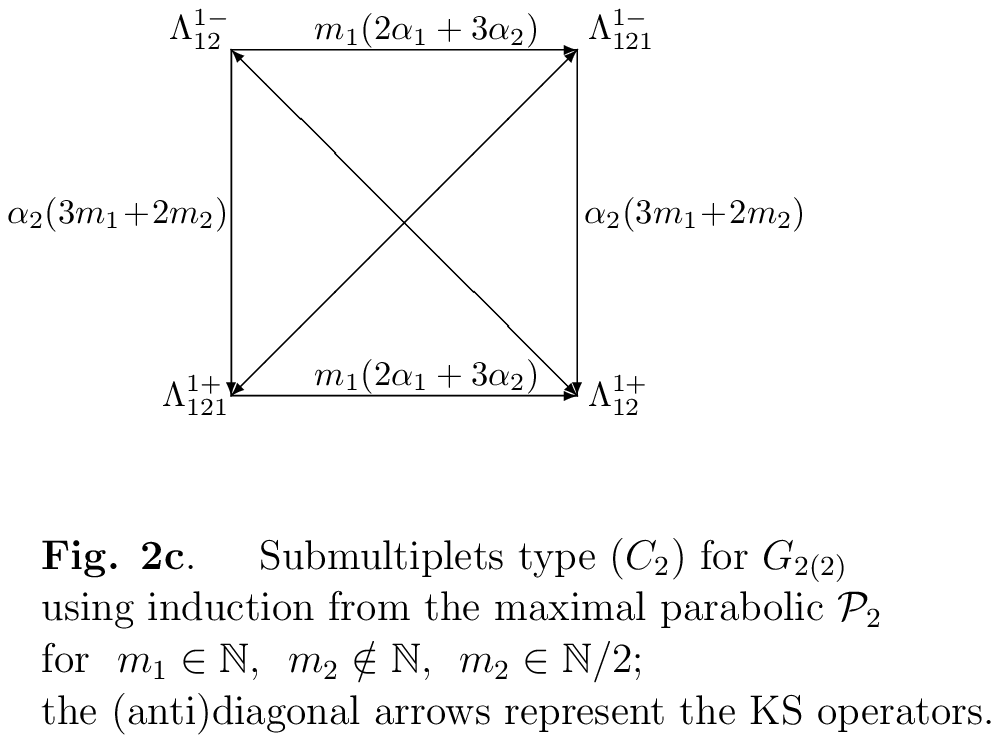}{20cm}

\end{document}